\documentclass[reqno, 11pt]{article}
\usepackage{amsthm}
\usepackage{amsmath}
\usepackage{amssymb}
\usepackage{amsfonts}
\usepackage[margin=1.5in]{geometry}
\usepackage{tikz}
\usepackage{tikz-cd}
\usetikzlibrary{arrows}
\usepackage{hyperref}
\usepackage{float}
\usepackage{subcaption}
\usepackage{wrapfig}
\newtheorem{definition}{Definition}
\newtheorem{theorem}{Theorem}
\newtheorem{corollary}{Corollary}
\newtheorem{proposition}{Proposition}

\newtheorem{lemma}{Lemma}
\newtheorem{remark}{Remark}
\newtheorem*{question}{Question}

\DeclareMathOperator{\Mod}{Mod}

\usepackage{comment}
\usepackage[shortlabels]{enumitem}
\newcommand{\N}{{\Bbb N}}
\newcommand{\C}{{\Bbb C}}

\newcommand{\qf}{{\mathcal Q}_f}
\newcommand{\q}{{\mathbf q}}
\newcommand{\cf}{{\mathcal C}_f}

\newcommand{\pf}{{\mathcal P}_f}

\title{Independence of Thurston eigenvalues and Lyapunov exponents for quadratic rational maps}
\author{Ayesha Bennett and Holly Krieger}
\date{\today}

\begin{document}

\maketitle
\begin{abstract}
    We study post-critically finite quadratic rational maps $f: \mathbb{P}^1\to\mathbb{P}^1$ where the two critical points lie in the same cycle. The \emph{Thurston eigenvalues} of $f$ are those of the derivative of the Thurston pullback map $\sigma_f$ at its fixed point. Addressing questions of Buff--Epstein--Koch, we construct a sequence of real, bitransitive centres in the critical orbit relation curve $C_2$ which degenerates in moduli space, and whose Lyapunov exponents tend to infinity whilst the Thurston eigenvalues satisfy $\frac{1}{10}< |\lambda|<1$. We also  show the Thurston eigenvalues across all bitransitive centres in $C_2$ have a spectral gap on compact subsets of the curve.
\end{abstract}

\section{Introduction}\label{Sintro}

Let $f\colon\mathbb{P}^1\to\mathbb{P}^1$ be a rational map of degree at least two, with finite post-critical set $P_f$. If $|P_f|\ge 4$ and $f$ is not a Latt\`es map, Thurston's pullback map $\sigma_f$ is a contracting holomorphic self-map of the Teichm\"uller space $\mathcal{T}_f$ of the sphere marked at the points of $P_f$, and it has a unique fixed point $\star$ \cite{doudhub}. We call the eigenvalues of the derivative $d_\star\sigma_f$ the \emph{Thurston eigenvalues} of $f$. Equivalently, these are the  eigenvalues of the pushforward map $f_*$ on $\mathcal{Q}_f$. We defer the related background to Section~\ref{SbackgroundC2}.

Buff, Epstein, and Koch studied the Thurston eigenvalues of unicritical polynomials and proved the existence of a spectral gap: every Thurston eigenvalue $\lambda$ of  $z^d+c$ with periodic critical point satisfies $\frac{1}{4d}<|\lambda|$ \cite[Theorem~1.4]{bek}. Their lower bound comes from the bound $|f'|\le 2d$ on the filled Julia set, and it is therefore natural to compare the eigenvalues with the Lyapunov exponent
\[
L(f)=\int_{\mathbb{P}^1}\log\|Df\|\,d\mu_f
\]
of $f$ with respect to its measure of maximal entropy $\mu_f$. For post-critically finite polynomials, $L(f)=\log d$ \cite[Section~11]{DeM03}, leading to the natural question \cite[Question~4]{bek}:

\begin{question}[Buff--Epstein--Koch]\label{q:BEK4}
For post-critically finite maps $f$, is there a relation between $\exp(-L(f))$ and $\inf_{\lambda}\{|\lambda|\}$, where $\lambda$ ranges over the Thurston eigenvalues of $f$?
\end{question}

In this paper, we show that the Lyapunov exponent can tend to infinity along a sequence while the Thurston eigenvalues stay bounded away from $0$. 

\begin{theorem}\label{cor:lyap}
There is no real constant $C>0$ such that $\inf_{\lambda}\{|\lambda|\}\le C\,e^{-L(f)}$ for all post-critically finite quadratic rational maps $f$.
\end{theorem}
In fact, Theorem~\ref{cor:lyap} is actually a corollary of the following result, which addresses another question raised in \cite[Question~3]{bek} about the relationship between bounded families of maps in moduli space and a spectral gap for the Thurston eigenvalues.
Denote by $\mathcal{M}_2$ the moduli space of quadratic rational maps: the quotient of the space of degree two rational maps $f: \mathbb{P}^1 \rightarrow \mathbb{P}^1$ by M\"obius conjugation. 
 
\begin{theorem}\label{maintheorem}
There exists a sequence $\{ f_{c_n}\}_{n \geq 2}$ of rational maps of degree 2 with periodic critical points such that
\begin{enumerate}
\item[(i)] the sequence $[f_{c_n}]$ leaves every compact subset of $\mathcal{M}_2$;
\item[(ii)] for each $n \geq 2$, every Thurston eigenvalue $\lambda$ of $f_{c_n}$ satisfies $\frac{1}{10}<|\lambda|<1$.
\end{enumerate}
\end{theorem}

The maps we consider are conjugate to the real post-critically finite quadratic rational maps of rotation number $(n+1)/(2n+1)$ constructed by Filom and Pilgrim \cite{filompilgrim}. Each $f_{c_n}$ is the centre of a hyperbolic component of bitransitive type in the sense of Milnor and Rees \cite{MilnorLei,Rees1990}: both critical points lie in the immediate basin of one attracting cycle, in different Fatou components.  We consider these maps as lying in the critical orbit relation curve $C_2$, the curve of maps for which one critical point maps to another after two applications of the map.  In coordinates, we choose
\[
f_c(z)=\frac{z^2-1}{z^2-c},
\]
for $c \in \mathbb{C} \setminus \{1\}.$ 

In Milnor's coordinates, the curve $C_2$ is a conic in $\mathbb{C}^2$ with two punctures, corresponding to $c=\infty$ and $c=1$. As $c\to1$, $[f_c]$ tends to the ideal point $\langle-1,-1,\infty\rangle$ of Milnor’s compactification of $\mathcal{M}_2$, at which two fixed-point multipliers tend to $-1$ and the third to $\infty$ \cite{MilnorLei}; our sequence converges to this ideal point. In contrast, post-critically finite unicritical polynomials are bounded in moduli space, since $c$ lies in the compact Multibrot set. In the quadratic case this can also be seen in Milnor’s coordinates: $z^2+c$ tends to $\langle0,\infty,\infty\rangle$ as $c\to\infty$, and post-critically finite maps can only accumulate at ideal points $\langle\mu,\mu^{-1},\infty\rangle$ with $\mu$ a root of unity.

The question of the relationship between a spectral gap of Thurston eigenvalues and boundedness in moduli space is really two questions, asking for bounds in either direction, and the bound in the other direction \emph{does} hold for the bitransitive centres on this curve $C_2$.

\begin{proposition}\label{prop:compactgap}
Let $K\subset\mathbb{C}\setminus\{0,1\}$ be compact. There is $\delta_K>0$ such that every Thurston eigenvalue $\lambda$ of every bitransitive centre $f_c$ with $c\in K$ satisfies $|\lambda|\ge\delta_K$.
\end{proposition}

This proposition is proved in Section \ref{cptsection}, and is an elementary consequence of the computation of the characteristic polynomial of $d_{\star}\sigma_f$ carried out in the proof of Theorem~\ref{maintheorem}. \\

{\bf Strategy of the proof.} The main theorem rests on a computation of the characteristic polynomial of $d_\star\sigma_f$, through a construction of Koch \cite{KOCH2013573}. The pullback map $\sigma_f$ does not descend to the moduli space $\operatorname{Mod}(P_f)$ of marked configurations of $P_f$, but, for the maps we consider, it covers a rational map $g_f$ in the opposite direction: $g_f\circ\pi\circ\sigma_f=\pi$ on $\mathcal{T}_f$, where $\pi\colon\mathcal{T}_f\to\operatorname{Mod}(P_f)$ is the universal covering. For a map $f=f_c$ in $C_2$ with periodic critical points, $g_f$ is given by a simple explicit formula in the coordinates of the critical cycle (Section~2.3). It is regular at $\pi(\star)$, and the chain rule gives
\[
Dg_f(\pi(\star))=d_\star\pi\circ(d_\star\sigma_f)^{-1}\circ(d_\star\pi)^{-1}.
\]
So the Thurston eigenvalues of $f$ are the reciprocals of the eigenvalues of the explicit matrix $Dg_f(\pi(\star))$ (in particular none of them is $0$), and the coefficients of its characteristic polynomial are, up to explicit elementary factors, products of the derivatives $f'(z_j)$ at points of the critical cycle (Section~4.1), in close analogy with \cite[Proposition~9.7]{bek}. The difficulty in proving Theorem~\ref{maintheorem} is the degeneration of the derivatives as $c\to1$; however, the \emph{products} of the derivatives along the critical cycle can potentially be controlled, since $f_c^2$ converges to the quadratic rational map
\[
g(z)=\frac{2}{3-z^2}
\]
locally uniformly away from $\pm1$ and $\pm\sqrt3$, a rescaling limit in the sense of Kiwi \cite{kiwi}; see also \cite{epstein,DeMarco}.
Rouch\'e's theorem then confines the roots of the characteristic polynomial of $Dg_f(\pi(\star))$ to the disc of radius $10$, which is Theorem~\ref{maintheorem}(ii).

The parameters $c_n$ come from the work of Filom and Pilgrim on real quadratic rational maps \cite{filompilgrim}: for each rational $p/q$ they construct a real critically finite quadratic rational map whose critical points lie in a cycle of period $q$ on $\widehat{\mathbb{R}}$, permuted as the rotation by $p/q$ \cite[Proposition~3.2]{filompilgrim}. To prove Theorem~\ref{maintheorem}(i),  we use the Pommerenke--Levin--Yoccoz inequality, in the form established by Petersen for a repelling fixed point on the boundary of a periodic Fatou component \cite[Theorem~D]{Petersen}. \\

{\bf Forthcoming work and further questions.} 
The goal of this paper is to present an interesting case of degeneration of Thurston eigenvalues independent from Lyapunov exponents, provable by elementary means.  In more generality, the degeneration $c \to 1$ is best understood through the language of trees of spheres \cite{Arf,DF}, or, equivalently, through the dynamics induced on the Berkovich projective line over the field of Puiseux series in $\varepsilon=1-c$ \cite{kiwi}. The tree has three vertices: the sphere $S_0$ centred at $0$ at scale $1$, and the spheres $S_{\pm1}$ centred at $\pm1$ at scale $|\varepsilon|$. The vertex $S_{-1}$ is strictly preperiodic, mapping to $S_0$, and the vertices $S_0,S_1$ form a two-cycle, with the first return map to $S_0$ being the rescaling limit $g$. The fixed point $z=1$ of $g$ is parabolic, and lies at the point of $S_0$ where $S_1$ is attached. The critical cycle of a bitransitive centre near $c=1$ moves through these vertices, and through the parabolic gate, along a constrained itinerary. To leading order in $\varepsilon$, $f'_c$ along the cycle are determined by this itinerary. So, therefore, are  the coefficients of the characteristic polynomial of Section~\ref{SScharpoly}, in which degenerating factors cancel in pairs along each return to $S_0$, reflecting the convergence $f_c^2 \to g$. The bitransitive centres of Theorem~\ref{maintheorem} were chosen as their itinerary is the simplest, and hence admit elementary study. The same tree controls the behaviour of the multipliers of the cycles of $f_c$ which pass through its vertices, and multipliers along degenerating sequences can be described by the non-archimedean techniques of Favre and Gong \cite{FG1,FG2} and of Luo \cite{Luo2,Luo1}. 

The forthcoming thesis work of the first author takes this point of view for a more general study of degeneration of Thurston eigenvalues for $C_2$ and other critical orbit curves in the space of quadratic rational maps \cite{Ben}.  Given this, while a number of interesting questions have been raised in the course of this work, we leave them out of this paper by choice. We will gladly discuss them with anyone who wishes to work on them in the collaborative spirit that the authors believe the subject deserves. \\

{\bf AI usage statement.} The ideas and mathematics in this paper were produced by the two authors. Figures were generated by OpenAI's ChatGPT-6.2 Astra.  A complete draft was given to Anthropic's Claude Fable 5.1 to check for notational and mathematical correctness, and superficial changes were made by the model.  The same model was used to check the literature for references, though any omissions are the responsibility of the authors, and make several choices about exposition. The model was also used to reconcile the spelling used by the authors, who (having learned to write in different hemispheres) struggled to agree on the spelling of centre. \\

{\bf Acknowledgements.} The initial ideas for this paper were inspired by lectures given by Sarah Koch in 2013 as part of an AMS Mathematics Research Communities programme, and the second author benefitted greatly from related discussions at the time with Sarah Koch, Adam Epstein, Curt McMullen, and Dylan Thurston. More recently, discussions with Caroline Davis, Laura DeMarco, Alex Kapiamba, and Jan Kiwi helped the authors make key decisions about the direction and scope of this paper. 
The first author was supported by the Additional Funding Programme for Mathematical Sciences, delivered by EPSRC (EP/V521917/1) and the Heilbronn Institute for Mathematical Research. The second author is partially supported by funds from the Leverhulme Trust. Finally, both authors would like to thank the Isaac Newton Institute for Mathematical Sciences, Cambridge, for support and hospitality during the programme Complex dynamics: interactions and influences, where work on this paper was undertaken. This work was supported by EPSRC grant EP/Z000580/1.

\section{The Thurston pullback and the curve $C_2$}\label{SbackgroundC2}

In this section we collect the background material and fix the notation used in the remainder of the paper: the Thurston pullback map and its relation to the pushforward operator on quadratic differentials, and the curve $C_2$ of quadratic rational maps in which our sequence of post-critically finite maps will be constructed.

\subsection{The Thurston pullback map}\label{SSthurston}

Let $S^2$ denote a topological two-sphere, equipped with an orientation, and fix an identification $S^2 = \mathbb{P}^1$.  For a finite set $P \subset S^2$ containing at least three points, the Teichm\"uller space $\mathcal{T}_P = \mathcal{T}(S^2, P)$ is the space of equivalence classes of orientation-preserving homeomorphisms $\phi: S^2 \rightarrow \mathbb{P}^1$, with $\phi_1$ equivalent to $\phi_2$ if there is a M\"obius transformation $\nu: \mathbb{P}^1 \rightarrow \mathbb{P}^1$ so that $\phi_1 = \nu \circ \phi_2$ on $P$, and $\phi_1$ is isotopic to $\nu \circ \phi_2$ relative to $P$.

The Teichm\"uller space admits a natural projection $\pi_P$ to the moduli space $\Mod(P)$ of injections of $P$ into $\mathbb{P}^1$ modulo post-composition by M\"obius transformations, defined by $\pi_P([\phi]) = [\phi|_P]$.  The projection $\pi_P$ is a holomorphic universal covering map, realising $\mathcal{T}_P$ as a complex manifold of dimension $|P| - 3$.  Note that the identification $S^2 = \mathbb{P}^1$ provides both $\mathcal{T}_P$ and $\Mod(P)$ with natural basepoints: the classes of the identity and of the inclusion, respectively.

It will be convenient to fix coordinates on $\Mod(P)$.  Write $r := |P| - 3$, and enumerate $P = \{p_0, p_1, p_2, q_1, \dots, q_r\}$.  Every class in $\Mod(P)$ has a unique representative $\iota: P \hookrightarrow \mathbb{P}^1$ satisfying $\iota(p_0) = \infty$, $\iota(p_1) = 1$, and $\iota(p_2) = 0$, and the assignment $[\iota] \mapsto (\iota(q_1), \dots, \iota(q_r))$ identifies $\Mod(P)$ with the open subset of $\mathbb{A}^r$ of points with pairwise distinct coordinates, none equal to $0$ or $1$.  Once an ordered triple $(p_0, p_1, p_2)$ has been fixed, we use this identification without further comment, regarding $\pi_P$ as a map $\mathcal{T}_P \rightarrow \mathbb{A}^r$.

Suppose now that $f: S^2 \rightarrow S^2$ is an orientation-preserving branched covering of degree at least two, with critical set $\cf$.  The \emph{post-critical set} of $f$ is
$$\pf := \bigcup_{n \geq 1} f^n(\cf),$$
and $f$ is \emph{post-critically finite} if $|\pf| < \infty$.  Assume from now on that $f$ is post-critically finite with $|\pf| \geq 3$, and write $\mathcal{T}_f := \mathcal{T}_{\pf}$.  By pulling back complex structures along $f$, we obtain the \emph{Thurston pullback map} $\sigma_f: \mathcal{T}_f \rightarrow \mathcal{T}_f$; concretely, given $\tau \in \mathcal{T}_f$, there is a rational map $F_{\tau}$, determined up to the choices of representatives, so that the following diagram commutes:
\begin{figure}
\centering
\[ \begin{tikzcd} 
  S^2 \arrow[r, "\sigma_f(\tau)"] \arrow[d, "f"]
    & \mathbb{P}^1 \arrow[d, "F_{\tau}"] \\
  S^2 \arrow[r, "\tau"]
&\mathbb{P}^1 \end{tikzcd} \]
\caption{Thurston pullback map}
\label{fig:tpm}
\end{figure}
Here and below we abuse notation, denoting a class in $\mathcal{T}_f$ and a representative homeomorphism by the same symbol.

\begin{definition} Post-critically finite branched covers $f$ and $g$ are \emph{combinatorially equivalent} if there is a pair of orientation-preserving homeomorphisms $\phi_1, \phi_2: (S^2, \pf) \rightarrow (S^2, \mathcal{P}_g)$ with $\phi_1$ isotopic to $\phi_2$ relative to $\pf$, so that the following diagram commutes:
\[ \begin{tikzcd}
  (S^2, \pf) \arrow[r, "\phi_2"] \arrow[d, "f"]
    & (S^2, \mathcal{P}_g) \arrow[d, "g"] \\
  (S^2, \pf) \arrow[r, "\phi_1"]
&(S^2, \mathcal{P}_g) \end{tikzcd} \]
\end{definition}

The Thurston pullback map is holomorphic and weakly contracting for the Teichm\"uller metric, contracting for a rational map with hyperbolic orbifold; in particular, when $|\mathcal{P}_f|\geq 5$.  In this case, the pullback map has a unique fixed point if and only if $f$ is combinatorially equivalent to a rational map (see \cite[Section~2]{doudhub}).  When $f$ is itself a post-critically finite rational map, the unique fixed point $\star$ of $\mathcal{T}_f$ has $F_{\star} = f$.

For a finite set $A \subset \mathbb{P}^1$ with $|A| \geq 3$, write $\mathcal{Q}(A)$ for the space of meromorphic quadratic differentials on $\mathbb{P}^1$ which are holomorphic on $\mathbb{P}^1 \setminus A$ and have at most simple poles at the points of $A$; $\mathcal{Q}(A)$ is a complex vector space of dimension $|A| - 3$.  We refer to \cite[Section~3]{bek} for background on meromorphic quadratic differentials and the pushforward operator in this context; in the notation of that paper, $\mathcal{Q}(A) = \mathcal{Q}(\widehat{\mathbb{C}}; A)$.  For a post-critically finite rational map $f$, define $\qf := \mathcal{Q}(\pf)$.  A non-constant rational map $F: \mathbb{P}^1 \rightarrow \mathbb{P}^1$ acts on meromorphic quadratic differentials $\q=q\, (dz)^2$ by the pushforward
$$F_* \q(y) := \sum_{g} g^* \q(y) = \sum_{x \in F^{-1}(y)}\frac{q(x)}{F'(x)^2}\, (dy)^2$$
where the sum ranges over the local branches $g$ of $F^{-1}$.  If $\q$ has at most simple poles, all contained in a set $A$, then $F_* \q$ has at most simple poles, contained in $F(A) \cup \mathcal{V}_F$, where $\mathcal{V}_F$ denotes the set of critical values of $F$ (see \cite[Section~3.3]{bek}).  In particular, for a post-critically finite rational map $f$, the pushforward restricts to a linear operator
$$f_*: \qf \rightarrow \qf.$$

The connection between $\sigma_f$ and $f_*$ is provided by the following (see \cite[Proposition~3.1 and 3.2]{doudhub}).

\begin{proposition}\label{DHcotangent} Let $\tau \in \mathcal{T}_f$ have representative $\phi$.  The cotangent space to $\mathcal{T}_f$ at $\tau$ is canonically isomorphic to $\mathcal{Q}(\phi(\pf))$, and under this identification the coderivative of $\sigma_f$ at $\tau$ is the pushforward operator $(F_{\tau})_*$, where $F_{\tau}$ as in figure~\ref{fig:tpm}.
\end{proposition}

Applying this at the fixed point of $\sigma_f$ for a rational map $f$, and recalling that a linear map and its transpose have the same characteristic polynomial, we obtain:

\begin{corollary}\label{computewithThurston} Let $f$ be a post-critically finite rational map and $\star$ the unique fixed point of $\sigma_f$.  Then the characteristic polynomial of $d_{\star}\sigma_f$ is the characteristic polynomial of $f_*: \qf \rightarrow \qf$.
\end{corollary}
This demonstrates that the spectrum of the derivative of the Thurston pullback map $\sigma_f$ at the unique fixed point coincides with the eigenvalue spectrum of $f_*$ on $\mathcal{Q}_f$. Observe that as $\mathcal{Q}(\mathcal{P}_f)$ is a finite-dimensional $f_*$-invariant subspace of meromorphic quadratic differentials on $\hat \C$,
one can also study the set of eigenvalues of the induced operator, which, in line with \cite{bek}, we denote $\Lambda_f$. For post-critically finite rational maps, \cite[Theorem~1.2]{bek} provides a description of $\Lambda_f$ in terms of the multipliers along cycles not contained in $\mathcal{P}_f$, and proves that $0 \in \Lambda_f$ for all $f$. Therefore, this spectrum is not related to $\sigma_f$, nor is the induced operator invertible, and so is not relevant to the Thurston eigenvalues that we are interested in. 
\begin{remark}\label{unitdiskremark} By the contraction principle, the eigenvalues of the pushforward acting on the space of all meromorphic quadratic differentials on $\mathbb{P}^1$ with only simple poles -- and hence on the invariant subspace $\qf$ -- have modulus at most $1$, with equality only if $f$ is a Latt\`es map with four post-critical points (see \cite[Lemma~1.1]{bek}).  The maps of interest in this paper are non-Latt\`es, so every eigenvalue of $f_*$ satisfies $|\lambda| < 1$, and the content of Theorem~\ref{maintheorem} is the complementary lower bound. 
\end{remark}

\subsection{The curve $C_2$}\label{SScurve}

 A quadratic rational map has three fixed points, counted with multiplicity, and the elementary symmetric functions $\sigma_1, \sigma_2, \sigma_3$ of the three fixed point multipliers are invariants of the conjugacy class. By the holomorphic fixed point formula they satisfy the relation $\sigma_3 = \sigma_1 - 2$.

\begin{theorem} \label{milnoriso}
\textup{(Milnor; \cite[Lemma~3.1]{MilnorLei})} The map $[f] \mapsto (\sigma_1(f), \sigma_2(f))$ defines an isomorphism $\mathcal{M}_2 \rightarrow \mathbb{A}^2$.
\end{theorem}

A quadratic rational map has exactly two critical points, and they are distinct, as are its two critical values.  For $n \geq 1$, we define $C_n \subset \mathcal{M}_2$ to be the locus of conjugacy classes of maps $f$ admitting a labelling $c_1, c_2$ of the critical points with
$$f^n(c_1) = c_2;$$
$C_n$ is an algebraic curve in $\mathcal{M}_2 \cong \mathbb{A}^2$.  We now describe the curve $C_2.$  By direct computation, one verifies the following.

\begin{lemma}\label{smoothc2} $C_2$ is the smooth conic
$$2\sigma_1^2 + \sigma_1 \sigma_2 + 4 \sigma_1 - 2 \sigma_2 + 16 =0$$
in $\mathbb{A}^2.$
\end{lemma}

We may thus parametrise $C_2$, and we do so as follows.  Suppose $[f] \in C_2$, with critical points labelled so that $f^2(c_1) = c_2$.  Then the three points $c_1$, $f(c_1)$, $c_2$ are pairwise distinct: if $f(c_1) = c_1$, then $f^2(c_1) = c_1 \neq c_2$, while if $f(c_1) = c_2$, then $c_2 = f^2(c_1) = f(c_2)$, so that both critical values of $f$ equal $c_2$.  Conjugating so that $(c_1, f(c_1), c_2) = (\infty, 1, 0)$, the map $f$ has $\cf=\{0,\infty\}$ and satisfies $f(\infty) = 1$ and $f(1) = f^2(\infty) = 0$, which forces
$$f(z) = f_c(z) := \frac{z^2 - 1}{z^2 - c}$$
for some $c \in \mathbb{C}$, with $f_c$ of degree two exactly when $c \neq 1$.  Each member of the family satisfies the dynamical portrait
$$\infty \xrightarrow{\ 2 \ } 1 \longrightarrow 0 \xrightarrow{\ 2 \ } \frac{1}{c} \longrightarrow \cdots,$$
with labels indicating a local degree greater than one.

Now and throughout the rest of the paper we write $f_c$ for the map above, and
$$\mathcal{F} := \left\{ f_c : c \in \mathbb{C} \setminus \{1\} \right\}$$
for the corresponding family.  $c \mapsto [f_c]$ is a closed embedding of $\mathcal{F}$ into $\mathbb{A}^2$ with image $C_2$.  We have 
$$\sigma_1([f_c]) = \frac{2(c+3)}{c-1}, \ \sigma_2([f_c]) = -\frac{4(c^2+c+2)}{c-1},$$
with punctures at $c = \infty$ and $c = 1.$

As $c \rightarrow \infty$, $\sigma_1$ and hence $\sigma_3 = \sigma_1 -2$ remains bounded but $\sigma_2$ does not, so some fixed point multiplier goes to infinity and another must go to $0$; that is, the puncture at $\infty$ realises the ideal point of multipliers $0, \infty, \infty.$ In particular, for $|c|$ sufficiently large, $f_c$ is not post-critically finite.

As $c \rightarrow 1$, the maps $f_c$ degenerate, with $\sigma_1, \sigma_2 \rightarrow \infty$ and $f_c$ converging to the constant map $1$ locally uniformly on $\mathbb{P}^1 \setminus \{ \pm 1 \}$.  As will be shown in Section~\ref{Sfp}, $c = 1$ is an accumulation point of post-critically finite parameters.

Recall that a rational map is \emph{hyperbolic} if the orbit of every critical point converges to an attracting cycle.  Hyperbolic maps form an open subset of $\mathcal{M}_2$, whose connected components are called \emph{hyperbolic components}.  Following Rees \cite{Rees1990} and Milnor \cite[Section~4]{MilnorLei}, the hyperbolic quadratic rational maps fall into four types: 
\begin{itemize}
    \item \emph{adjacent}, with both critical points in the same Fatou component;
    \item \emph{bitransitive}, with the critical points in different components of the immediate basin of a single attracting cycle; 
    \item \emph{capture}, with only one critical point in the immediate basin of the unique attracting cycle;
    \item \emph{disjoint}, with two distinct attracting cycles.
\end{itemize}

Adjacent, bitransitive, and capture parameters occur among the hyperbolic components meeting $C_2$, while disjoint type is obstructed by the critical orbit relation.

Suppose now that $c$ is a parameter for which $f_c \in \mathcal{F}$ is a bitransitive centre of a hyperbolic component, and let $p$ denote the period of the superattracting cycle.  Writing $z_j := f_c^j(\infty)$, the portrait above shows that the cycle is
$$\infty = z_0 \xrightarrow{\ 2 \ } 1 = z_1 \longrightarrow 0 = z_2 \xrightarrow{\ 2 \ } z_3 \longrightarrow \cdots \longrightarrow z_{p-1} \longrightarrow z_0,$$
so that $p \geq 4$, and the post-critical set of $f_c$ is precisely this cycle.  Thus $|\pf| = p$, and
$$\dim \mathcal{T}_{f_c} = p - 3 =: r.$$

For a bitransitive centre $f_c$, every cycle disjoint from the critical cycle is repelling; in particular all three fixed points of $f_c$ are repelling, and their multipliers enter into the Yoccoz-style inequality of Section~\ref{Ssyoccoz}, applied at a fixed point on the boundary of an immediate basin.

\subsection{The inverse map on moduli}\label{SSgf}

We fix once and for all the coordinates on $\Mod(\pf) \subset \mathbb{A}^r$ described in Section~\ref{SSthurston} determined by the ordered triple $(p_0, p_1, p_2) := (z_0, z_1, z_2) = (\infty, 1, 0)$ and by $q_j := z_{j+2}$ for $1 \leq j \leq r$, and write $x_1, \dots, x_r$ for the resulting coordinates.  Since $f_c$ is rational, the unique fixed point $\star$ of $\sigma_{f_c}$ is the basepoint of $\mathcal{T}_{f_c}$, and its image in $\mathbb{A}^r$ is the configuration of the marked points themselves:
$$\pi(\star) = (z_3, z_4, \dots, z_{p-1}) \in \mathbb{A}^r.$$
In particular, at the fixed point the coordinate $x_r = z_{p-1}$ satisfies $f_c(x_r) = \infty$, and so
\begin{equation}\label{xrsquared}
    x_r^2 = c.
\end{equation}

Fix a bitransitive centre $f = f_c$ in $\mathcal{F}$, with post-critical cycle of period $p = r + 3 \geq 4$.  The pullback map $\sigma_f$ does not descend to a self-map of moduli; it does, however, cover a rational map in the \emph{opposite} direction, a phenomenon studied systematically by Koch \cite{KOCH2013573}.  For the family $\mathcal{F}$, this map admits the following explicit description, which will be used in Section~\ref{Sspecgap}.

\begin{lemma}\label{gfstructure} Let $f=f_c $ be a bitransitive centre in $\mathcal{F}$ with post-critical cycle of period $p = r + 3$.  The rational map $g_f: \mathbb{A}^r \dashrightarrow \mathbb{A}^r$ defined by
$$g_f(x_1, \dots, x_r) := \left( \frac{1}{x_r^2},\ f_{x_r^2}(x_1),\ \dots,\ f_{x_r^2}(x_{r-1}) \right) = \left( \frac{1}{x_r^2},\ \frac{x_1^2 - 1}{x_1^2 - x_r^2},\ \dots,\ \frac{x_{r-1}^2 - 1}{x_{r-1}^2 - x_r^2} \right)$$
is the unique rational map satisfying
$$g_f \circ \pi \circ \sigma_f = \pi$$
on $\mathcal{T}_f$; in particular, $g_f$ is defined at every point of $\pi(\sigma_f(\mathcal{T}_f))$.  Moreover, $g_f$ is regular at $\pi(\star)$, with $g_f(\pi(\star)) = \pi(\star)$.
\end{lemma}

\begin{proof} Let $\tau \in \mathcal{T}_f$, write $\tau' := \sigma_f(\tau)$, and choose representatives $\phi$ of $\tau$ and $\psi$ of $\tau'$, normalised as in Section~\ref{SSthurston} so that both send $(z_0, z_1, z_2) = (\infty, 1, 0)$ to $(\infty, 1, 0)$; thus
$$u := \pi(\tau) = (\phi(z_3), \dots, \phi(z_{p-1})) \quad \text{and} \quad v := \pi(\tau') = (\psi(z_3), \dots, \psi(z_{p-1})),$$
and we may choose the rational map $F_{\tau}$ of degree two so that $\phi \circ f = F_{\tau} \circ \psi$.  Evaluating this relation at the post-critical cycle gives
\begin{equation}\label{eqmarkedorbit}
    F_{\tau}(\psi(z_j)) = \phi(z_{j+1}), \qquad j \in \mathbb{Z}/p\mathbb{Z} \, .
\end{equation}
The critical points of $F_{\tau}$ are $\psi(\infty) = \infty$ and $\psi(0) = 0$, and \eqref{eqmarkedorbit} with $j = 0$ and $j = 1$ gives $F_{\tau}(\infty) = 1$ and $F_{\tau}(1) = 0$.  Arguing exactly as in Section~\ref{SScurve}, these conditions force $F_{\tau} = f_b$ for a unique $b \in \mathbb{C}$, with $b \neq 1$ since $F_{\tau}$ has degree two.

We now read off the remaining instances of \eqref{eqmarkedorbit}.  Taking $j = p - 1$: since $\phi(z_0) = \infty$, the coordinate $v_r = \psi(z_{p-1})$ satisfies $f_b(v_r) = \infty$, so
$$v_r^2 = b.$$
Taking $j = 2$: $u_1 = F_{\tau}(0) = f_b(0) = \dfrac{1}{b} = \dfrac{1}{v_r^2}$.  Taking $j = k + 2$ for $1 \leq k \leq r - 1$:
$$u_{k+1} = F_{\tau}(\psi(z_{k+2})) = f_{v_r^2}(v_k).$$
Thus $u = g_f(v)$, provided the right-hand side is defined at $v$.  It is: $v_r \neq \psi(z_2) = 0$ as $\psi$ is injective on $\pf$, so the first component is defined. Suppose that some later denominator vanished, say $v_k^2 = v_r^2$ with $k \leq r - 1$, then, since $f_b(-w) = f_b(w)$ ,we would have $\phi(z_{k+3}) = F_{\tau}(v_k) = F_{\tau}(v_r) = \infty = \phi(z_0)$, contradicting the injectivity of $\phi$ on $\pf$, as $z_{k+3} \neq z_0$ for $k \leq r - 1$.  As $\tau$ was arbitrary, $g_f \circ \pi \circ \sigma_f = \pi$ on $\mathcal{T}_f$.

Finally, $\sigma_f(\star) = \star$, so $\pi(\star) \in \pi(\sigma_f(\mathcal{T}_f))$; by the above, $g_f$ is defined there, and $g_f(\pi(\star)) = g_f(\pi(\sigma_f(\star))) = \pi(\star)$.

To prove uniqueness, simply note that since $\pi$ is a covering map, the identity $g_f \circ \pi \circ \sigma_f = \pi$ implies that $\pi(\sigma_f(\mathcal{T}_f))$ is Zariski dense in $\mathbb{A}^r$.

\end{proof}

Given the preceding, we now apply the chain rule to the identity $g_f \circ \pi \circ \sigma_f = \pi$ at $\star$. This gives
$$Dg_f(\pi(\star)) \circ d_{\star}(\pi \circ \sigma_f) = d_{\star}\pi.$$
Since $\pi$ is a covering map, $d_{\star}\pi$ is invertible, with
$$Dg_f(\pi(\star)) = (d_{\star}\pi) \circ (d_{\star}\sigma_f)^{-1} \circ (d_{\star}\pi)^{-1}.$$

Combined with Corollary~\ref{computewithThurston}, the eigenvalues of $f_*: \qf \rightarrow \qf$ are therefore precisely the reciprocals of the eigenvalues of the Jacobian matrix $Dg_f(\pi(\star))$; in particular, $f_*$ is invertible for every bitransitive centre in $\mathcal{F}$.  The resulting computation of the characteristic polynomial of $f_*$ is carried out in Section~\ref{Sspecgap}.

The construction of the bitransitive centres relevant to Theorem~\ref{maintheorem}, and the geometry of their critical cycles, is the subject of the next section.

\section{Geometry of Real, Bitransitive centres}\label{Sfp}
In this section, we explicitly describe the sequence of bitransitive, post-critically finite maps within $\mathcal{F}$ which are the subject of Theorem~\ref{maintheorem}.
\begin{figure}[H]
\centering
\begin{subfigure}[c]{0.39\textwidth}
  \centering
  \includegraphics[height=3.9in]{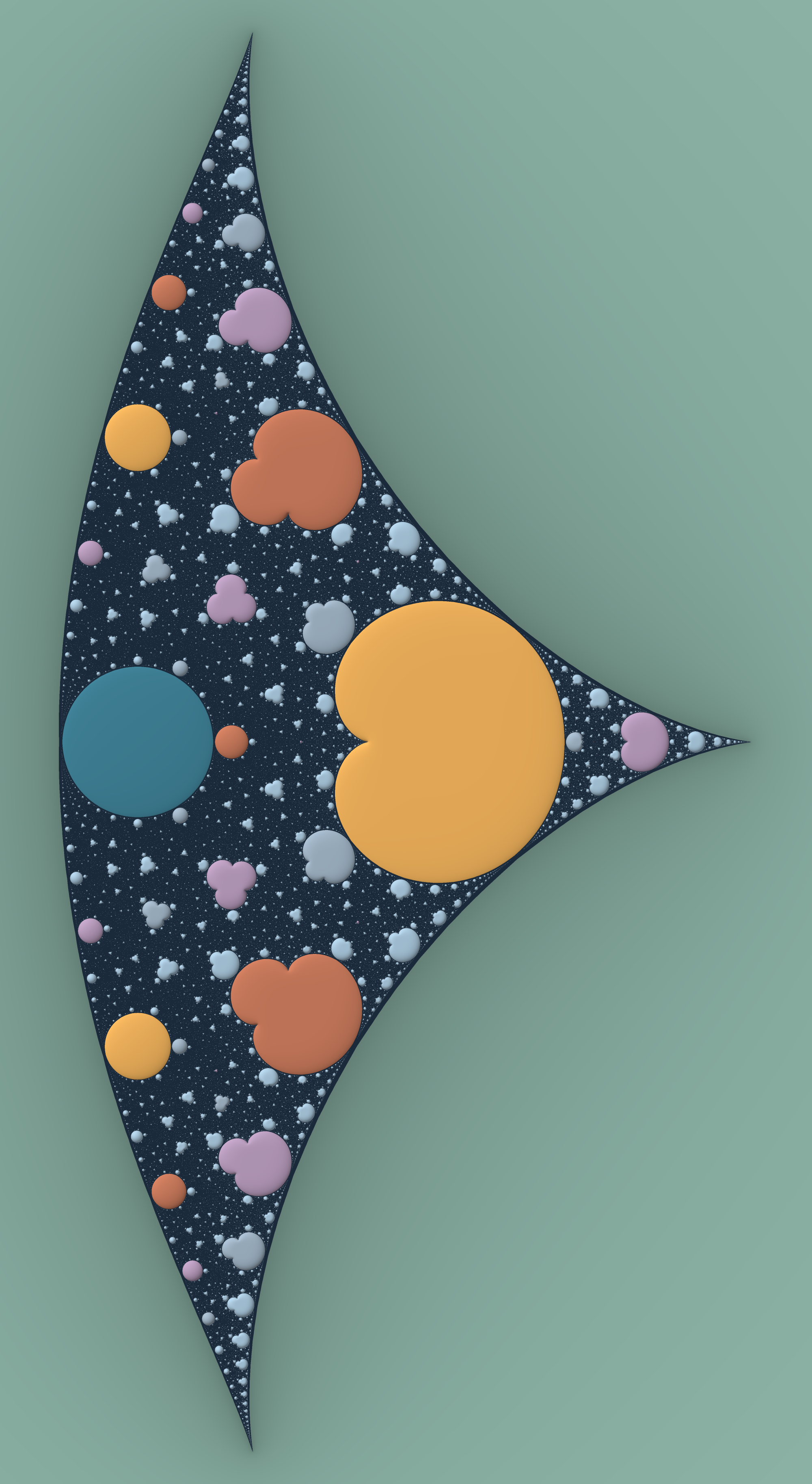}
  \caption{The parameter plane of $f_c$.}
  \label{fig:plane}
\end{subfigure}\hfill
\begin{minipage}[c]{0.58\textwidth}
  \begin{subfigure}{\linewidth}
    \centering
    \includegraphics[width=\linewidth]{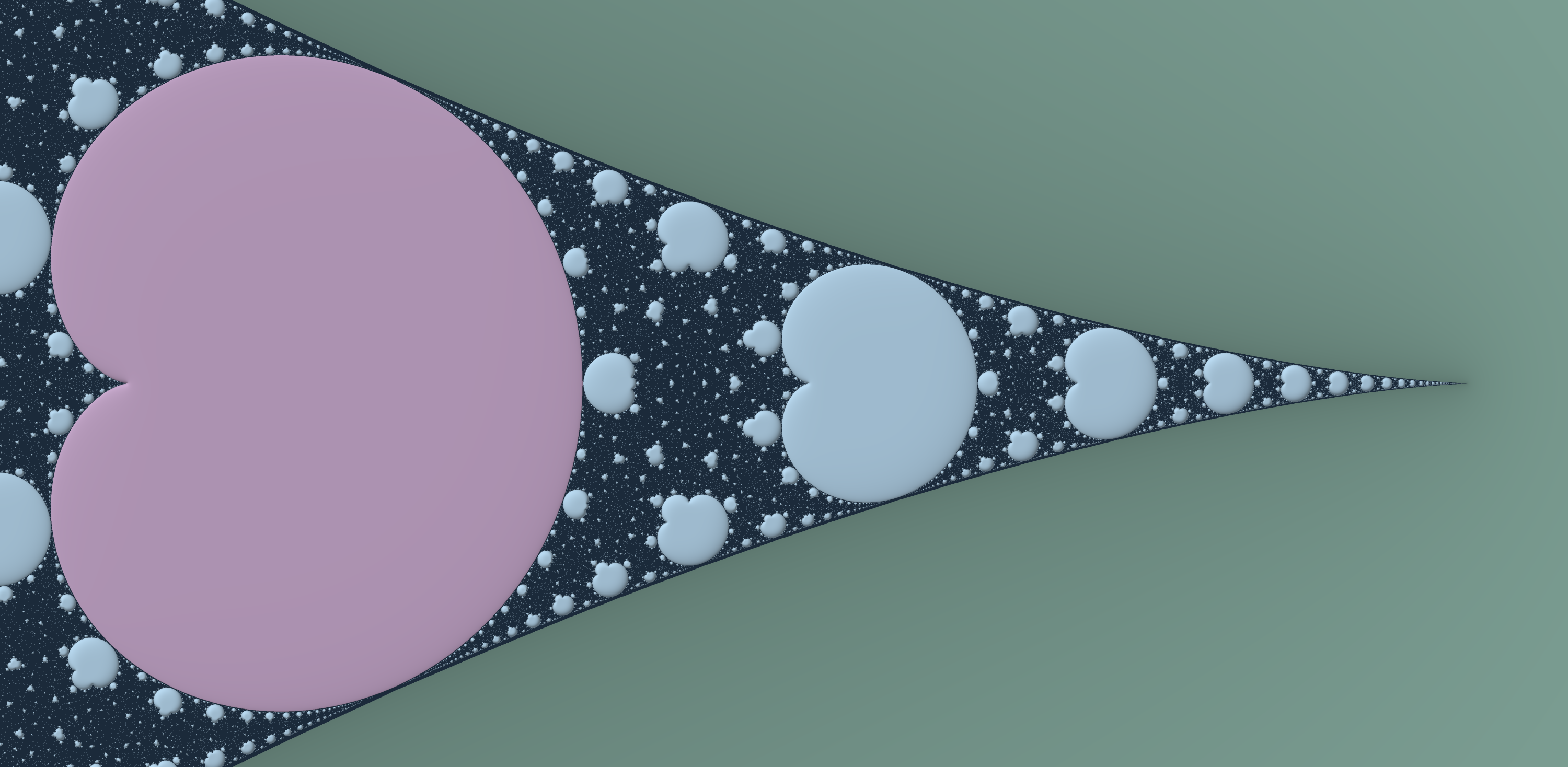}
    \caption{Zoom in near the puncture $c=1$.}
    \label{fig:zoom}
  \end{subfigure}\\[1ex]
  \begin{subfigure}{\linewidth}
    \centering
    \includegraphics[width=\linewidth]{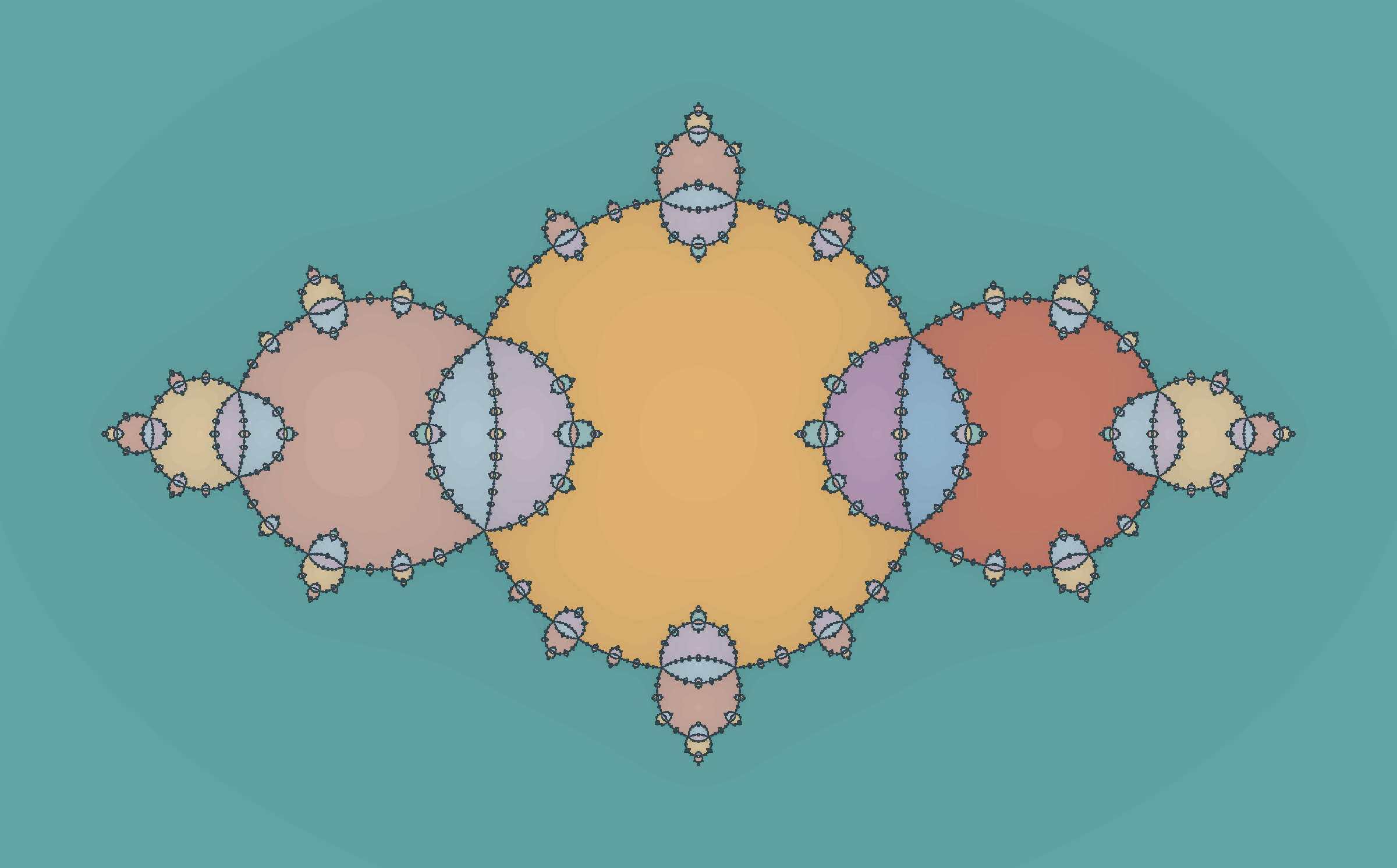}
    \caption{The Julia set of the period-$5$ centre $f_{c_2}$.}
    \label{fig:julia}
  \end{subfigure}
\end{minipage}
\caption{\small Parameter plane of $f_c$ and the Julia set for $c=c_2$.}
\label{fig:C2}
\end{figure}
In (a) and (b), the hyperbolic components are coloured according to the period of the attracting cycle capturing the critical points: period $2$ teal, $3$ orange, $4$ red, $5$ purple, and higher periods light blue. In (c), the critical points $0$ and $\infty$ lie in the orange and teal Fatou components respectively, and the cycle of immediate basins is orange $\to$ red $\to$ purple $\to$ teal $\to$ blue $\to$ orange, containing $0\mapsto 1/c_2\mapsto\sqrt{c_2}\mapsto\infty\mapsto1$. 
\subsection{Existence of real bitransitive centres}\label{Sscnreal}
The goal of this section is to establish the following proposition: 
\begin{proposition}\label{realcn}
For every $n \geq 2$, there exists $c_n \in (0,1)$ such that $f_{c_n}=\frac{z^2-1}{z^2-c_n}$ is bitransitive, post-critically finite map with critical points that lie in a period $2n+1$ cycle with rotation number $\frac{n+1}{2n+1}$. Furthermore, the critical portrait has the following ordering on $\hat{\mathbb{R}}:=\mathbb{R}\cup \{\infty\}$:
\begin{equation}\label{ordering}
0=c_2<z_4<\dots <z_{2n}<1<z_3<\dots<z_{2n-1}<\infty=c_1 \, .
\end{equation}
\end{proposition}

To prove the existence of these special parameters, we shall use \cite[Proposition~3.2]{filompilgrim}, which we state in full below for the convenience of the reader. 

\begin{proposition}\label{the family}
Let $q\geq 2$ be an integer and $p/q\in\mathbb{Q}/\mathbb{Z}$ a fraction in lowest terms. Let $p'$ denote the multiplicative inverse of $p$ mod $q$.

There exists a unique post-critically finite quadratic rational map $r$ (depending on $\frac{p}{q}$) with the following properties:
\begin{enumerate}[(a)]
\item \label{a} (marked critical points) the two critical points $\hat{c}_0, \hat{c}_1$ are labelled;
\item \label{b} (the real condition) $r$ has real coefficients;
\item \label{c} (normalisation) $\hat{c}_0=0$, $r({\rm{i}})={\rm{i}}, r(-{\rm{i}})=-{\rm{i}}$;
\item \label{d} (bitransitive) the critical points $\hat{c}_0, \hat{c}_1$ lie in a cycle of period $q$, so that $\mathcal{P}_r$ is given by $\left\{r^{ j}(\hat{c}_0)\right\}_{j=0}^{q-1}$;
\item \label{e} (rotation number) $\mathcal{P}_r $ is a subset of circle $\hat{\Bbb{R}}$ equipped with its usual orientation; alternatively,  we may write its elements  as $\{\hat{x}_0, \hat{x}_1, \ldots, \hat{x}_{q-1}\}$ where 
\begin{equation}\label{auxiliary 2}
 0 =\hat{c}_0=\hat{x}_0< \hat{x}_1< \hat{x}_2 < \ldots < \hat{x}_{q-1}<0;
\end{equation}
then we have 
$$  r(\hat{x}_j)=\hat{x}_{j+p}$$ 
for all $j \in \Bbb{Z}/q\Bbb{Z}$, so that the restriction $r\restriction_{\mathcal{P}_r}: \mathcal{P}_r \rightarrow \mathcal{P}_r$ forms a cycle with rotation number $p/q$;
\item \label{f} (adjacent critical points) $\hat{c}_1=\hat{x}_1$, which in conjunction with the previous property implies 
$r^{ p'}(\hat{c}_0)=\hat{c}_1$ \,.

\end{enumerate}
Furthermore, the map satisfies two additional properties
\begin{enumerate}[(I)]
    \item \label{i} (basins and their boundaries) for each $j\in\Bbb{Z}/q\Bbb{Z}$, the immediate super-attracting basin $\Omega_j$ containing $\hat{x}_j$ is a Jordan domain, and
$
\partial\Omega_j\supset{\pm{\rm i}} \, .
$
\item \label{ii} (normalised basin coordinates) there exist unique Riemann maps (extended to the boundaries)
$\phi_j:(\overline{\Bbb{D}}, 0, 1)\rightarrow(\overline{\Omega_j}, \hat{x}_j, \zeta_j), j \in \Bbb{Z}/q\Bbb{Z}$.  The first-return map of the immediate basin $\Omega_0$ to itself,  in these coordinates,  is 
$$\phi_j^{-1}\circ r^{q}\circ\phi_j: (\Bbb{D}, 0, 1)\rightarrow(\Bbb{D}, 0, 1);$$
 and is given by $\omega\mapsto\omega^4$.
\end{enumerate}
\end{proposition}
We shall refer to the maps $r$ satisfying the above properties, as Filom-Pilgrim maps. Notice that, although these maps satisfy most of the properties stated in Proposition~\ref{realcn}, they have a different normalisation to those in our family $\mathcal{F}$ (see Section~\ref{SScurve}). The proof of Proposition~\ref{realcn} will demonstrate the relationship between $r$ (for an appropriate choices of $p,q$) and certain real, positive bitransitive centres of $\mathcal{F}$. 
 
\begin{proof}[Proof of Proposition~\ref{realcn}] Let $r_{p/q}$ denote the Filom-Pilgrim map with rotation $p/q$, and consider the case $p/q = \frac{n+1}{2n+1}$ for $n \geq 2$.  In this case the multiplicative inverse $p'$ of $p \mod q$ is $p' = 2,$ so by property \ref{f} of Proposition \ref{the family}, $r_{(n+1)/(2n+1)}^2(0)$ is a critical point and the map is conjugate to an element of $C_2$.  

Let $\hat{x}_i$ be as in Proposition \ref{the family} property \ref{e} and let $\mu$ be the unique M\"obius map which takes $(\hat{x}_0=0, \hat{x}_1, \hat{x}_p)$ to $( \infty, 0, 1)$.  Then $\mu \circ r_{(n+1)/(2n+1)} \circ \mu^{-1} \in C_2$, and we write it as $f_{c_n}$ for some $c_n \in \mathbb{C} \setminus \{ 1 \}$.  Since $\mu$ is orientation-preserving, the inequality \ref{ordering} follow.

\end{proof}

\subsection{Convergence of real bitransitive centres}\label{Ssyoccoz}

We continue to write $c_n$ for the real, bitransitive centres defined in Proposition~\ref{realcn}. From here, and the rest of the paper, we shall write $f_n$ to denote $f_{c_n}$. We now wish to show that they converge to a puncture in $\mathcal{F}$; in particular:
\begin{proposition}\label{cnto1} The sequence $c_n$ converges to $1$ as $n \rightarrow \infty.$
\end{proposition}
Note that this is expected from the parameter space -- see figure~\ref{fig:zoom}

The result follows almost immediately once we have established that some non-real fixed point has multiplier tending to $-1$ as $n \to \infty$. Taking this as fact, the proof is as follows:
\begin{proof}[Proof of Proposition~\ref{cnto1}] Let $\beta_n:=\beta(c_n)$ be the marked fixed point of $f_n$ associated to the fixed point $i$ of the Filom-Pilgrim map. Note that this means $\beta_n$ is not real.  By Proposition~\ref{cnmulsto1}, the finite multipliers $f_{n}'(\beta_n) \rightarrow -1$ as $n \rightarrow \infty$, and the same is true for the fixed point $\overline{\beta_n}$ since $f_{n}$ has real coefficients.  Suppose that $\{ c_n \}_{n \in \mathbb{N}}$ has an accumulation point in $[0, 1)$, and call it $c_{\infty}.$  Then $f_{{\infty}}$ has two fixed points of multiplier $-1$, violating the holomorphic Lefschetz fixed point formula. Hence, the sequence must converge to $1$ as claimed.
\end{proof}

Therefore, it remains to understand the behaviour of the multipliers as the sequence degenerates. We shall use \cite[Theorem~D]{Petersen} stated below, which may be seen as a generalisation of the Yoccoz inequality for certain rational maps.
\begin{theorem}\label{Peterson}
  Let $\Lambda$ be a simply connected $m$-periodic, hyperbolic or parabolic Fatou domain associated to a rational map $R$. Then there exists a constant $M_\Lambda$, depending only on the restriction of  $R^m:\Lambda \to \Lambda$ such that for any  accessible boundary point $\alpha \in \partial \Lambda$ with period $k \in \N$, multiplier $\lambda \in \C \setminus \mathbb{D}$ and rotation number $\frac{p}{q}$, there exists a logarithm $L \in \mathbb{H}_+$ of $\lambda$ with 
  \begin{equation}\label{petest}
      \left|L- 2\pi i\frac{p}{q}\right| \leq \frac{2k\sin(\theta)}{qmN} M_\Lambda \, ,
  \end{equation}
  where $\theta$ is the angle between $2\pi i$ and $L-2\pi i\frac{p}{q} $ and $N$ is the number of cycles of periodic access to $\alpha$ represented in $\Lambda$. Moreover,
 if $\Lambda$ is hyperbolic. Then $M_\Lambda \leq \log(d)+D$ where $d=\deg(R^m:\Lambda \to \Lambda)$ and $D$ is the maximum hyperbolic distance in $\Lambda$ between the attracting fixed point for $R^m$ in $\Lambda$ and its pre-images in $\Lambda$. 
\end{theorem} 
With this in hand, we can prove the assumed `fact' used in the proof of Proposition~\ref{cnto1}, thus concluding the section:

\begin{proposition}\label{cnmulsto1}
The map $f_n$ has a non-real fixed point $\beta_n$ with multiplier tending to $-1$ as $n \to \infty$.
\end{proposition}
\begin{proof}
As $f_n$ is conjugate to the Filom-Pilgrim map $r_{(n+1)/(2n+1)}$ by a real M\"obius map, it suffices to prove this statement for $r_{(n+1)/(2n+1)}.  $
    Now, $i$ is fixed by $r_{(n+1)/(2n+1)}$ and by Proposition~\ref{the family}, property~\ref{i} and property~\ref{ii}, there is a cycle of accesses at $i$, which the map rotates with rotation number $\frac{n+1}{2n+1}.$ Applying Theorem~\ref{Peterson} with $d=4$ and $D=0$ (as follows from Proposition~\ref{the family} property~\ref{ii}), taking $\Lambda$ to be the immediate basin $\Omega_0$ of $0$, we have from \eqref{petest} that:
      \begin{eqnarray*}
                \left|L- 2\pi i\frac{p}{q}\right| &\leq& \frac{2k\sin(\theta)}{q^2N} M_\Lambda \\ & \leq &
                \frac{2}{(2n+1)^2}\log(4) \, .
      \end{eqnarray*}
This goes to $0$ as $n \to \infty$, and so $L \to i\pi$, and so we conclude $\lambda \to -1$. 

\end{proof}

\section{Existence of Spectral Gap}\label{Sspecgap}

In this section, we first indirectly compute the characteristic polynomial of the pushforward map $(f_c)_*: \mathcal{Q}_{f_c} \rightarrow \mathcal{Q}_{f_c}$ via study of a related rational map arising from the inverse of $\sigma_f$. We then establish upper bounds for the coefficients of this map, enabling us to estimate the eigenvalues of $(f_{c})_*$ and prove the existence of a spectral gap. 

\subsection{Computing the Characteristic Polynomial}\label{SScharpoly}
The characteristic polynomial of $f_*$, which we denote by $\xi_f$ can be computed directly with respect to a basis of $\mathcal{Q}_f$. This is the approach taken in \cite{bek}, where the authors prove a spectral gap for the pushforward operator associated to a family of unicritical polynomials with a periodic critical point. In this paper, we take the approach as in Koch's thesis \cite{kochthesis}, discussed in Section~\ref{SSgf}, which explicitly exploits properties of bitransitive type maps, and bypassing working explicitly with quadratic differentials. 
As demonstrated in Corollary~\ref{computewithThurston}, we can study the Thurston pullback map $\sigma_f$ on Teichm\"uller space to deduce properties of $f_*$.  Furthermore, as discussed in Section~\ref{SSgf}, for bitransitive, post-critically finite quadratic maps, there is a unique rational map $g_f: \mathbb{A}^r \dashrightarrow \mathbb{A}^r$ arising from the inverse of $\sigma_f$, which is regular at the $\pi-$image of the unique fixed point $\star$ in the Teichm\"uller space. By Lemma~\ref{gfstructure} it suffices to compute the characteristic polynomial of $Dg_f(\pi(\star))$.  Note that at $\pi(\star)$, $x_r$ satisfies \eqref{xrsquared}, with $x_r \in (0,1)$ for $c=c_n$, and so $f_n=f_{x_r^2}$.

\begin{lemma}\label{charpoly}
    Let $r \geq 2$ and $f=f_c$ for any bitransitive centre $c$.  Then the characteristic polynomial of $Dg_f(\pi(\star))$, denoted by $p$, is given by 
   \begin{align*}
p(\lambda)={}&\lambda^r
-\frac{2x_r(x_{r-1}^2-1)}{(x_{r-1}^2-x_r^2)^2}\lambda^{r-1}
-\sum_{j=1}^{r-2}
\frac{2x_r(x_j^2-1)}{(x_j^2-x_r^2)^2}
\left(\prod_{i=j+1}^{r-1} f'(x_i)\right)\lambda^j
+\frac{2}{x_r^3}\prod_{i=1}^{r-1}f'(x_i).
\end{align*}
    Furthermore, it satisfies the following relation with $\xi_f$: \begin{equation}\label{relcharpolys}
        p(0)\xi_f(\lambda)=\lambda^r p(1/\lambda) \, .
    \end{equation}
\end{lemma}
\begin{proof}
Let $\mathbf{x}^* = (x_1, \dots, x_r) = \pi(\star)$ denote the unique fixed point image of $\sigma_f$. Recall from Lemma~\ref{gfstructure} that the components of $g_f(x_1, \dots, x_r) = (y_1, \dots, y_r)$ are given by
$$
y_1 = \frac{1}{x_r^2}, \qquad y_j = \frac{x_{j-1}^2 - 1}{x_{j-1}^2 - x_r^2} \quad \text{for } j = 2, \dots, r.
$$
We first compute the entries of the $(r \times r)$ Jacobian matrix $J = Dg_f(\mathbf{x}^*)$. Differentiating each component yields:
\\

\noindent For $y_1$:
    $$
    \frac{\partial y_1}{\partial x_r} = -\frac{2}{x_r^3}, \qquad \frac{\partial y_1}{\partial x_j} = 0 \quad \text{for } 1 \le j \le r-1.
    $$  
 For $y_j$ ($2 \le j \le r$):
    \begin{align*}
    \frac{\partial y_j}{\partial x_{j-1}} &= \frac{2x_{j-1}(1 - x_r^2)}{(x_{j-1}^2 - x_r^2)^2} = f_c'(x_{j-1}), \\[1ex]
    \frac{\partial y_j}{\partial x_r} &= \frac{2x_r(x_{j-1}^2 - 1)}{(x_{j-1}^2 - x_r^2)^2}, \\[1ex]
    \frac{\partial y_j}{\partial x_k} &= 0 \quad \text{for } k \notin \{j-1, r\}.
    \end{align*}

Thus, the matrix $(\lambda I-J)$ has the following form:
$$
\lambda I - J = \begin{pmatrix}
\lambda & 0 & 0 & \dots & 0 & \frac{2}{x_r^3} \\[1ex]
-a_1 & \lambda & 0 & \dots & 0 & -b_2 \\[1ex]
0 & -a_2 & \lambda & \dots & 0 & -b_3 \\
\vdots & \vdots & \ddots & \ddots & \vdots & \vdots \\
0 & 0 & \dots & -a_{r-2} & \lambda & -b_{r-1} \\[1ex]
0 & 0 & \dots & 0 & -a_{r-1} & \lambda - b_r
\end{pmatrix}
$$
where $a_j := \frac{\partial y_{j+1}}{\partial x_j} = f'(x_j)$ for $1 \le j \le r-1$, and $b_j := \frac{\partial y_j}{\partial x_r} = \frac{2x_r(x_{j-1}^2 - 1)}{(x_{j-1}^2 - x_r^2)^2}$ for $2 \le j \le r$. We want to compute the characteristic polynomial $p(\lambda) = \det(\lambda I - J)$.

There are many ways to do this, and we obtain
$$
p(\lambda) = (\lambda - b_r)\lambda^{r-1} - \sum_{j=1}^{r-2} b_{j+1} \left( \prod_{i=j+1}^{r-1} a_i \right) \lambda^j +  \frac{2}{x_r^3}\prod_{j=1}^{r-1} a_j.
$$
Substituting the expressions for $b_{j+1} $ and $a_j$ yields our desired form.

Finally, the relation $p(0) \xi_f(\lambda) = \lambda^r p(1/\lambda)$ follows directly from the relationship $g_f, \sigma_f$ and $f_*$ from Corollary~\ref{computewithThurston}.
\end{proof}

For notational simplicity, we will write  \begin{equation}\label{charpolyshort}
    p(\lambda) = \sum_{j = 0}^r \alpha_j \lambda^j \, ,
\end{equation} 
where $\alpha_j$ are the coefficients of $p(\lambda)$ explicitly written in Lemma~\ref{charpoly}. 

Our next goal is to obtain upper bounds on $|\alpha_j|$. Once these are established, we can apply Rouch\'e's Theorem on a sufficiently large circle, with a choice of functions that will imply that all $r$ roots of $p(\lambda)$ lie inside a disc of uniformly bounded radius. Using \eqref{relcharpolys}, we can then translate this estimate into a spectral gap for $\xi_f$, and hence obtain a lower bound on the modulus of the eigenvalues of $f_*$ 
\\

To bound the coefficients of the characteristic polynomial $p(\lambda)$, we need to control products of derivatives along the critical cycle. Although individual derivatives $f'(x_i)$ may become arbitrarily large as the parameter approaches the puncture, we can use the theory of rescaling limits to help establish bounds on the products.

\subsection{Rescaling limits}\label{SSrl}

Rescaling limits were introduced by Stimson in his thesis \cite{stim} and developed extensively in complex dynamics by Epstein \cite{epstein}, Kiwi \cite{kiwi}, and DeMarco \cite{DeMarco}. For our quadratic rational family $f_c(z) = \frac{z^2-1}{z^2-c}$, the map collapses to the constant function $1$ as $c \to 1$. However, the second iterate $f_c^2$ admits a non-trivial algebraic limit.

\begin{lemma}\label{lem:rescaling_limit}
As $c \to 1$, the second iterate sequence $f_c^2(z)$ converges uniformly on compact subsets of $\mathbb{C} \setminus \{\pm 1, \pm \sqrt{3}\}$ to the quadratic rational map
\begin{equation}\label{rlmap}
    g(z) := \lim_{c \to 1} f_c^2(z) = \frac{2}{3-z^2}.
\end{equation}
\end{lemma}

\begin{proof}
Direct evaluation of the second iterate yields
$$
f_c^2(z) = \frac{(\frac{z^2-1}{z^2-c})^2 - 1}{(\frac{z^2-1}{z^2-c})^2 - c} = \frac{(z^2-1)^2 - (z^2-c)^2}{(z^2-1)^2 - c(z^2-c)^2}.
$$
Expanding the terms in both numerator and denominator and factoring out $(1-c)$ gives
$$
f_c^2(z) = \frac{(1-c)(-2z^2+c+1)}{(1-c)(z^4 -2(c+1)z^2+c^2+c+1)} = \frac{-2z^2+c+1}{z^4-2(c+1)z^2+c^2+c+1}.
$$
Taking $c \to 1$, the expression simplifies  to
$$
g(z) = \frac{-2z^2 + 2}{z^4 -4z^2 +3} = \frac{-2(z^2-1)}{(z^2-1)(z^2-3)} = \frac{2}{3-z^2}.
$$
Uniform convergence on compact subsets avoiding the poles of $g$ and the cancelled zeros/poles at $\pm 1$ follows.
\end{proof}

The critical portrait of the limiting map $g$ is illustrated in figure~\ref{fig:critical_portrait}. 

\begin{figure}[H]
\centering
\begin{tikzpicture}[>=stealth, node distance=1.8cm, auto]
  \node (p1) {$\infty$};
  \node (p2) [right of=p1] {$0$};
  \node (p3) [right of=p2] {$\frac{2}{3}$};
  \node (dots) [right of=p3] {$\cdots$};
  \node 
  (one) [right of=dots] {$1$};

  \draw[->] (p1) -- node[above] {\scriptsize $2$} (p2);
  \draw[->] (p2) -- node[above] {\scriptsize $2$} (p3);
  \draw[->] (p3) -- (dots);
  \draw[->] (one) to[out=60,in=120,looseness=6] (one);
\end{tikzpicture}
\caption{Critical portrait of the second iterate rescaling limit $g(z) = \frac{2}{3-z^2}$.}
\label{fig:critical_portrait}
\end{figure}
Note that the orbit of $g^n(\infty)$ approaches the non-degenerate parabolic fixed point $1$ as $n \to \infty$.
The unique pole at $z = \sqrt{3}$ on the positive real axis corresponds precisely to the limit of $x_{r-1}$ as $x_r \to 1$. Consequently, for all orbit points $x_i$ bounded strictly away from $x_{r-1}$, the derivative $(f^2)'(x_i)$ remains uniformly bounded by $g'(x_i)$.

\subsection{Bounding the Eigenvalues}\label{SSboundeig}
In this section, we prove the existence of a uniform bound on the coefficients $|\alpha_j|$ from \eqref{charpolyshort}; that is, there exists $R>0$ so that $|\alpha_j| < R^{r-j}$ for each $0 \leq j \leq r-1$.  We will continue to abuse notation and write $f$ for $f_n=f_{x_r^2}$.
\\

Observe that all coefficients of $p(\lambda)$ involve products of the derivative along the critical cycle, except $\alpha_r$ and $\alpha_{r-1}$. We start by bounding these two. Clearly, $\alpha_r=1$. It is straightforward to bound $|\alpha_{r-1}|$:

\begin{lemma}\label{ceasy} $|\alpha_{r-1}| < 1 \, .$
\end{lemma}

\begin{proof} Since $\frac{x_{r-1}^2-1}{x_{r-1}^2-x_r^2} = x_r,$ we have $x_{r-1}^2 = 1+x_r+x_r^2.$  Therefore 
$$-\alpha_{r-1} = \frac{2x_r(x_r+x_r^2)}{(1+x_r)^2} = \frac{2x_r^2}{1+x_r}.$$
Since $x_r \in (0, 1)$, we have $|\alpha_{r-1}| < 1$ as desired.
\end{proof}

We turn to the remaining coefficients $|\alpha_j|$, $ j \in \{0,\dots,r-2\}$. First, we use the rescaling limit \eqref{rlmap} discussed in Section~\ref{SSrl} to bound the product of derivatives along the critical orbit. The heuristic is as follows: the convergence to the rescaling limit map $g$, given in \eqref{rlmap}, for real values of $z$ is fairly rapid away from $x_r < z <1$.   Thus, $(f^{2})'(z)$ behaves similarly to $g'(z) = \frac{4z}{(3-z^2)^2}$ (which is bounded on intervals not containing $\pm \sqrt{3}$) away from this interval. Therefore, outside the bad interval, we can utilise the rescaling limits for any $x_i$ not converging to the poles of $g$ as the parameter $x_r$ approaches the puncture at $1$. 
 Recall from Proposition~\ref{realcn} that $x_i \in (0,\infty)$ for all $i$, and so they avoid the interval $(x_r,1)$. Furthermore, they obey a certain ordering-- \eqref{ordering}.

The following lemma makes this heuristic precise.
\begin{lemma}\label{rlbound}
    Suppose that $x_i$ is in the critical orbit of $f$. Then for all  $i \neq r-1$, we have $|(f^2)'(x_i)| \leq 12 \sqrt{3}$.  
\end{lemma}

\begin{proof}
As we are working at the fixed point  of $\sigma_{f}$, we have $x_r^2=c_n$. Therefore, 
$$f^2(x) = \frac{x_r^2+1-2x^2}{(x^2-x_r^2+x_r-1)(x^2-x_r^2-x_r-1)} \, .$$
Computing the derivative, we find
$$h(x) := (f^2)'(x) = \frac{4x(x^2-1)(x^2-x_r^2)}{(x^2-x_r^2+x_r-1)^2(x^2-x_r^2-x_r-1)^2}.$$
We want upper bounds on $h(x_i)$. For $ x \in \mathbb{R}_{\geq0}$, $h(x)$ has two asymptotes: $$0<\sqrt{x_r^2-x_r+1}<1 \: \text{and} \: 1<\sqrt{x_r^2+x_r+1}=x_{r-1}  \, ,$$ where the last equality comes from the relationship $x_r=\frac{x_{r-1}^2-1}{x_{r-1}^2-x_r^2}$. The fact that $x_{r-1}$ is a pole for $h(x)$ is precisely why it is excluded from the statement of Lemma~\ref{rlbound}. Moreover, $x_{r-1}^2=1+x_r+x_r^2<3$, so $x_{r-1}<\sqrt3$. Dynamically, this reflects the fact that as $x_r\to1$, $x_{r-1}$ increases to $\sqrt3$, the pole of the rescaling limit $g$.
\\

We now bound $h(x_i)$ for $i\neq r-1$ for $x_i \in (0,1)$. In this interval, $h$ is positive from $(0,x_r)$. We will show that 
$$\max\limits_{x_i \in (0,x_r)}h(x_i)<1 \, .$$ To compute this maxima, we first rewrite $h(x)$ in a more useful form
\begin{equation}\label{hgen}
    h(x) =  4x(x^2-1)(x^2-x_r^2)    \frac{f^2(x)^2}{(x_r^2+1-2x^2)^2} \, .
\end{equation}
By the Arithmetic Mean - Geometric Mean Inequality, noting the denominator can be expressed as $(-(x^2-1)-(x^2-x_r^2))^2$,  we can bound $h(x_i)$:
\begin{eqnarray*}
    h(x) &\leq&  x \, f^2(x)^2 \frac{4(1-x^2)(x_r^2-x^2)}{((1-x^2)+(x_r^2-x^2))^2}
    \\ &\leq& x f^2(x)^2 \, .
\end{eqnarray*}
Finally, rewriting
\begin{eqnarray*}
f^2(x)&=& \frac{-2x^2+x_r^2+1}{x^4-2(1+x^2_r)x^2+x_r^4+x_r^2+1} \\[1ex]&=& \frac{-2x^2+x_r^2+1}{-2x^2+x_r^2+1+(x^2-x^2_r)^2} \, ,
\end{eqnarray*}
we can easily see that the denominator is strictly greater than the numerator for $0<x\neq x_r$, and both are positive in the interval of interest. Therefore $h(x)<1$ for $x \in (0,x_r)$. By considering the expression for $h$ given in \eqref{hgen}, it is clear that for the same $x$, $h(x)>0$. Hence, we have $0<h(x)<x<x_r<1$, and so $0<h(x_i)<1$ for $x_i \in (0,x_r)$.
This concludes finding an upper bound on $h(x_i)$ for $x_i \in (0,1)$.
\\

We now move to establishing an upper bound for $h(x_i)$ for $x_i \in (1,\infty)$. First, note that by \eqref{ordering},  we have $x_i \leq x_{r-1} \leq \sqrt{3}$ for all $i$. 

As on $(1, x_{r-1})$,  $h(x)$ is positive and increasing, the maximum value along the critical orbit (excluding at $x_{r-1}$, the asymptote) is at $x_{r-3}$. Working with $h(x_{r-3})$ in the form of \eqref{hgen}, we can factor $(x^2_{r-3}-x_{r}^2)^2$ out of the denominator and use the relationship $x_{r-2} = \frac{x_{r-3}^2-1}{x_{r-3}^2-x_r^2} \, $ twice, to obtain:
$$
h(x_{r-3})=\frac{4x_{r-3}x_{r-2}x^2_{r-1}}{(x_{r-2}+1)^2} \, .
$$
Hence, we have $h(x_{r-3}) < 4x^3_{r-1}< 12 \sqrt{3}$ as desired.

\end{proof}
 \begin{remark}
     Within this proof, we have shown that $|h(x_i)| <1$ for even cycle elements, i.e all $x_i \in (0,1)$.  Upon observing that for  $ x \in (1,x_{r-1})$, the function $f^2_{x_r}(x)$ is increasing and greater than $1$, we also have that $1<|h(x_i)|< 12\sqrt{3} $ for odd cycle elements, i.e all $x_i \in (1,\infty)$, excluding $x_{r-1}$.
 \end{remark}

With Lemma~\ref{rlbound} in hand, we can bound the remaining coefficients of $p(\lambda)$.
\begin{corollary}\label{cj} Let $1 \leq j \leq r-2$.  Then $|\alpha_j| <5^{r-j}$.
\end{corollary}

\begin{proof} First, suppose that $j$ is even. Then, by taking out the term $f'(x_{r-1})$ from the product, we have
$$\alpha_j =  \frac{2x_r(x_j^2-1)}{(x_j^2-x_r^2)^2} \prod_{i=j+1}^{r-1} f'(x_i) = \frac{4x_{r-1}x_r(1-x_r^2)(x_j^2-1)}{(x_j^2 -x_r^2)^2(x_{r-1}^2-x_r^2)^2} \prod_{i=j+1}^{r-2} f'(x_i).$$
 There now are an even number of terms in the product, and so by the chain rule and Lemma~\ref{rlbound}, 
\begin{eqnarray*}
|\alpha_j| &<& \left( \frac{4x_{r-1}x_r(1-x_r^2)(x_j^2-1)}{(x_j^2 -x_r^2)^2(x_{r-1}^2-x_r^2)^2} \right) (12\sqrt{3})^{(r-j-2)/2} 
\\[1ex] & \leq& \left(\frac{4x_{r-1}x_r}{(x_{r-1}^2-x_r^2)^2}\right) x_{j+1}  \left(\frac{1-x_r^2}{x_r^2-x_j^2}\right) (12\sqrt{3})^{(r-j-2)/2} \\[1ex] &\leq& \left(\frac{4x_{r-1}x_r}{(1+x_r)^2}\right) x_{j+1}  (1-x_{j+1}) (12\sqrt{3})^{(r-j-2)/2}  \, ,
\end{eqnarray*}
where, for the last inequality, we have used the fact that
$$\frac{1-x_r^2}{x_j^2-x_r^2} = 1 - x_{j+1} \: \, \text{and} \: \,\frac{x_{r-1}^2-1}{x_{r-1}^2-x_{r}^2} = x_r \, .$$
Therefore, as $x_r < 1 < x_{j+1} \leq x_{r-1} < \sqrt{3}$, we have
$$|\alpha_j| < 12 (x_{j+1}-1)(12\sqrt{3})^{(r-j-2)/2} < (12\sqrt{3})^{(r-j)/2} < 5^{r-j},$$
as claimed. \\

Suppose now that $j$ is odd.  We rearrange so that the product again begins at an odd index:
$$\alpha_j =  \frac{2x_r(x_j^2-1)}{(x_j^2-x_r^2)^2} \prod_{i=j+1}^{r-1} f'(x_i) = \frac{x_r(x_j^2-1)}{x_j(1-x_r^2)} \prod_{i=j}^{r-1} f'(x_i).$$
Again we pull out $f'(x_{r-1})$ to compensate for the term external to the product, obtaining
$$\alpha_j = \frac{2x_{r-1}x_r(x_j^2-1)}{x_j(x_{r-1}^2-x_r^2)^2}  \prod_{i=j}^{r-2} f'(x_i).$$
To bound $|\alpha_j|$ from above, we use the following observations. First, as $j$ is odd, $x_r < 1 < x_j < x_{r-1} < \sqrt{3}$. Second, we use $x_{r-1}^2-x_r^2=1+x_r>1$ to bound the $(x_{r-1}^2-x_r^2)^2$ term in the denominator. Finally, we use that $x-\frac{1}{x}$ is an increasing function for $x>0$, and so, upon combining these observations with Lemma~\ref{rlbound}, we have
$$|\alpha_j| < 4 (12\sqrt{3})^{(r-j-1)/2} <  (12\sqrt{3})^{(r-j)/2}<5^{r-j} \, .$$
\end{proof}
We bound the remaining constant coefficient.
\begin{corollary}\label{c0}
   The coefficient $\alpha_0$ satisfies the upper bound $|\alpha_0|<5^r$.
\end{corollary}
\begin{proof}
First, using $$x_{r-1}^2=1+x_r+x_r^2 \, ,$$ and our ordering  \eqref{ordering}, we have $$\frac{1}{x_r^2}=x_1 \leq x_{r-1}$$
and so 
$$
1 \leq x_r^4(1+x_r+x_r^2) \, 
$$
with $x_r^4(1+x_r+x_r^2)$ increasing on $(0,1)$.
If we put in $x_r=0.8$, then we get $1 \leq 0.9994$,  a contradiction, and so for $c=x_r^2$ to be a real bitransitive centre, we must have $x_r>0.8$.
\\
We use Lemma~\ref{rlbound}, after pulling the factor $f'(x_{r-1})$ out of the product as before:
    $$|\alpha_0| \leq \frac{2}{x_r^3} \left( \frac{2x_{r-1}(1-x_r^2)}{(x_{r-1}^2-x_r^2)^2} \right)(12\sqrt{3})^{(r-2)/2} < \frac{2}{0.8^3} (2 \sqrt{3}) 5^{r-2} < 5^r \, .$$

\end{proof}
Finally, using the method set out after Lemma~\ref{charpoly}, we shall prove Theorem~\ref{maintheorem}.

\begin{proof}[Proof of Theorem~\ref{maintheorem}]
    First note that by Proposition \ref{cnto1}, the sequence $f_n$ leaves any compact subset of the quadratic moduli space, satisfying the first item of Theorem \ref{maintheorem}.  So we must prove the eigenvalue bound, and will show that all the roots of the characteristic polynomial $p(\lambda)$ from Lemma~\ref{charpoly} are bounded above independently of $r$. To do this, we use Rouch\'e's Theorem with the functions $\lambda^r- \alpha_{r-1}\lambda^{r-1}$ and $p(\lambda)-\lambda^r+\alpha_{r-1}\lambda^{r-1}$. By the bounds of Lemma~\ref{ceasy}, Corollary~\ref{c0} and Corollary~\ref{cj}, $|\alpha_j| < 5^{r-j}$ for $0 \leq j \leq r-1$, so by Rouch\'e's theorem, all eigenvalues have modulus strictly less than 10.

By \eqref{relcharpolys}, we deduce the lower bound of Theorem~\ref{maintheorem}.
\end{proof}

We finish the section by proving Theorem \ref{cor:lyap}.  The result follows from a combination of Theorem~\ref{maintheorem}, and an application of result \cite[Corollary~1.6]{DeM03}, stated below for convenience, to our setting, which provides an explicit formula for the Lyapunov exponent:
\begin{lemma}\label{laura}
    The Lyapunov exponent of a rational map f of degree $d > 1$
is given by
$$
L(f)=\sum_{n=1}^{2d-2}G_F(c_j)-\log(d)-\frac{2}{d}\log{|\text{\emph{Res}}(F)|}\, ,
$$
where $F : \C^2 \to \C^2$
is any homogeneous polynomial map such that $\pi \circ F =
f \circ \pi$, $G_F=\lim_{n \to \infty} \frac{\log||F^n(z)||}{\deg(F)^n}$ is the escape rate function of $F$ and points $c_j \in \C^2$
are determined by the condition $| \det DF(z)| =
\prod_j|z \wedge c_j |$
\end{lemma}
By applying Lemma~\ref{laura} in our setting, it is straightforward to obtain a lower bound on the Lyapunov exponent along our sequence of maps. Indeed, combining the result of the below corollary  with Theorem \ref{maintheorem} immediately proves Theorem~\ref{cor:lyap}.
\begin{corollary}\label{cor:lyapeff}
For every $n\ge2$,
\[
L(f_{c_n})\ \ge\ \tfrac18\log\frac{1}{|1-c_n|}-\frac{1}{4}\log2 .
\]
In particular $L(f_{c_n})\to\infty$ along the sequence of Theorem~\ref{maintheorem}.
\end{corollary}

\begin{proof}
Given $n \geq 2,$ lift $f_{c_n}$ to $$F(z,w)=(F_1,F_2)=(z^2-w^2,\,z^2-c_nw^2) \, ,$$ put the sup norm on $\mathbb{C}^2$, so then $G_F(v)=\lim_{k\rightarrow \infty} 2^{-k}\log\|F^k(v)\|$. 
Since $\det DF= 4wz(1-c)$, we may take $c_1=(0,1)$ and $c_2=(4(1-c),0)$. As Res($f$)$= (1-c)^2$, and $G_F(\lambda v)= G_F(v)+\log|\lambda|$, by Lemma~\ref{laura}, we have
\begin{eqnarray*}
L(f_{c_n})&=&G_F(0,1)+G_F(1,0)+\log2-\log|1-c_n| \\&=&5\,G_F(1,0)+\log2-2\log|1-c_n| \, ,
\end{eqnarray*}
where we have used $F^{ 2}(1,0)=(0,1-c_n)$ for the second equality.

To bound $G_F(1,0)$ from below, note that since $c_n \in (0,1)$, the identities $F_1-F_2=(c_n-1)w^2$ and $c_nF_1-F_2=(c_n-1)z^2$ show that $\|F(z,w)\|\ \ge \frac{|1-c_n|}{4}$ for every unit vector $(z,w)$. Writing $u_k=F^k((1,0))/\|F^k((1,0))\|$, we have
$$G_F((1,0))=\frac{1}{2} \log 1 + \frac{1}{4} \log |1-c_n| + \frac{1}{8} \log 1+ \sum_{k\ge3}\frac{\log\|F(u_k)\|}{2^{k+1}},$$
since $\log \|F(u_k)\| = \log \|F^{k+1}((1,0))\| - 2 \log \| F^k((1,0)) \|$ by homogeneity.  Thus
$$G_F((1,0))\geq \frac{1}{4} \log |1 - c_n| +( \log |1-c_n| - \log 4)  \sum_{k\ge3}2^{-k-1} = \frac{3}{8} \log |1 - c_n| - \frac{1}{4} \log 2,$$
and so
$$L(f_{c_n}) \geq \frac{1}{8} \log \frac{1}{|1-c_n|} - \frac{1}{4} \log 2$$
as claimed.
\end{proof}

\section{A bound in the compact case} \label{cptsection}

In this section we prove Proposition \ref{prop:compactgap} of the introduction:

\begin{proof}[Proof of Proposition~\ref{prop:compactgap}]
We use the notation of $x_i$ for the critical orbit and the characteristic polynomial as given in Lemma~\ref{charpoly}: 
\begin{align*}
p(\lambda)={}&\lambda^r
-\frac{2x_r(x_{r-1}^2-1)}{(x_{r-1}^2-x_r^2)^2}\lambda^{r-1}
-\sum_{j=1}^{r-2}
\frac{2x_r(x_j^2-1)}{(x_j^2-x_r^2)^2}
\left(\prod_{i=j+1}^{r-1} f'(x_i)\right)\lambda^j
+\frac{2}{x_r^3}\prod_{i=1}^{r-1}f'(x_i).
\end{align*}
We will show that if $c=x_r^2$ lies in a compact subset of the curve, then the modulus of the coefficients satisfy an upper bound of the form $C^{r-j}$  and hence, by Rouch\'e's Theorem, obtain a bound on the eigenvalues.

Assume throughout that $\delta \leq |c-1| \leq R$ for some positive $\delta, R$, and say a complex number $x$ is \emph{large} if $|x| \geq M$, where we will choose $M$ later.

We will use the following identities numerous times:
\begin{equation}\label{usefuleq}
    x_{i+1} - 1 = \frac{x_r^2 - 1}{x_i^2 - x_r^2} \, ,
\end{equation}
and also \begin{equation}\label{usefulderiveq}
    f_c'(x_i) = \frac{2 x_i (1-x_r^2)}{(x_i^2-x_r^2)^2} = \frac{2x_i(x_{i+1}-1)^2}{1-c} \, ,
\end{equation}
where we use \eqref{usefuleq} for the second equality. 
\\

If $x_i$ is large, then using \eqref{usefuleq}, it is easy to see that $|f_c(x_i)|<2$ for any  $M^2\geq 2(R+1)$, and hence the image is not large. In particular, no two consecutive points in the critical cycle can be large.
Then using the first equality from \eqref{usefulderiveq} we have $|f_c'(x_i)| \ll 1/|x_i|^3.$

Now for $x_{i}$ not large, since $|c-1| \geq \delta$, using the second equality from \eqref{usefulderiveq}, we have that $$|f_c'(x_i)| \ll |x_i|| x_{i+1}-1|^2 \quad \text{and} \quad \ll |x_i||x_{i+1}|^2
\,\quad \text{for $x_{i+1}$ large}\,    .$$
We move to bounding products of derivatives $\prod_{i = j}^{r-1} |f_c'(x_i)|.$ Let $I$ be the set of indices in this range for which $x_{i+1}$ is large.  Then $I$ has no two successive indicies as noted above, and since $|x_r^2| = |c|<M^2$, $x_r$ is not large, so $I \subset \{ j, \dots, r-2 \}.$ 
\\
For each $i \in I$ we have that 
\begin{eqnarray}\label{xinotlarge}
x_i \text{ is not large} &\implies&  |f_c'(x_i)| \ll |x_{i+1}|^2 \, \\ \text{and also } \  x_{i+1} \, \text{is large} &\implies& \label{xi+1large}  |f_c'(x_{i+1})| \ll \frac{1}{|x_{i+1}|^3} \, .
\end{eqnarray}
So, in the product of derivatives, the pair of terms with indices $i, i+1$ with $i \in I$ have product which is, by \eqref{xi+1large},\eqref{xinotlarge},  $\ll 1/|x_{i+1}|$ and hence uniformly bounded.  

For $i \not\in I$, we have a uniform upper bound on $|f_c'(x_i)|$ as follows: if $x_i$ is large then $|f_c'(x_i)| \ll 1/|x_i|^3 = O(1)$, while if $x_i$ is not large, then neither $x_i$ nor $x_{i+1}$ is large, and so $|f'_c(x_{i})| \ll |x_i||x_{i+1}-1|^2\leq M(M+1)^2=O(1)$.  Thus we have some uniform constant $C$ such that 
$$\prod_{i = j}^{r-1} |f_c'(x_i)| \leq C^{r-j}$$
for all $1 \leq j \leq r.$

Using compactness of $K \subset \C \setminus \{0,1\}$, and the fact it does not contain $0$, the modulus of the constant coefficient $|\alpha_0|$ from Lemma~\ref{charpoly} is then bounded above by $C^{r-1}$ for some constant $C$, and for the remaining coefficients $1 \leq j \leq r - 1$, since
$$\alpha_j = -\frac{2x_r(x_j^2-1)}{(x_j^2-x_r^2)^2} \left(\prod_{i=j+1}^{r-1} f'_c(x_i)\right),$$
we have bounds as follows. If $j = r-1$ then the product is the empty product $1$.  If $x_{j + 1}$ is not large, then $x_j$ is not close to $\pm x_r$. Then if $x_j$ is not large, our bound on the product coefficient is immediate, whilst if $x_j$ is large then we use the bound   $\gg \frac{1}{2}|x_j|^2$ for the term in the denominator. Therefore, for $x_{j+1}$ not large, the term $-\frac{2x_r(x_j^2-1)}{(x_j^2-x_r^2)^2}$ is uniformly bounded.    
For the case that $x_{j+1}$ is large,  note that as
\begin{equation}\label{usefuleq2}
    x_j^2 - x_r^2 = \frac{x_j^2-1}{x_{j+1}}
\end{equation}
we have
$$x_r^2 - 1 = x_j^2 - 1 - \frac{x_j^2 - 1}{x_{j+1}} = (x_j^2 - 1)(1 - \frac{1}{x_{j+1}})\, ,$$
and so for $x_{j+1}$ large, $|x_j^2 - 1|$ is uniformly comparable to $|x_r^2 - 1| \geq \delta$.  Using \eqref{usefuleq2},
$$-\frac{2x_r(x_j^2-1)}{(x_j^2-x_r^2)^2} =-\frac{2x_rx^2_{j+1}}{x_j^2-1}$$ 
and so, upon using \eqref{xi+1large},
$$ \left| \frac{2x_r(x_j^2-1)}{(x_j^2-x_r^2)^2}  f_c'(x_{j+1}) \right| \ll \frac{|x_{j+1}|^2}{|x_{j+1}|^3}$$ is uniformly bounded as well.  Thus, in either case, we have a bound of the form $C^{r - j}$ on $\alpha_j$ for all $1 \leq j \leq r-1$.  From Rouch\'e's Theorem, the bound on eigenvalues follows.

\end{proof}

\newpage


\end{document}